\documentclass[amstex,thmsa,  twocolumn]{article} 
\usepackage[tbtags]{amsmath}
\usepackage[]{amssymb}
\usepackage[]{amsbsy}
\usepackage[]{amsthm}
\usepackage{graphicx}

\newtheorem{theorem}{Theorem}[section]

\begin{document}


\title{\Large \textbf{The exponential tracking and disturbance rejection for the unstable Burgers' equation
with   general references and disturbances}}

\author{    Weijiu Liu\thanks{Corresponding author. Email:
weijiul@uca.edu, Phone: 1-501-450-5661,
 Fax: 1-501-450-5662}      \\
\normalsize{ Department of  Mathematics }\\
 \normalsize{  University of Central Arkansas }\\
  \normalsize{201 Donaghey Avenue, Conway, AR 72035, USA }\\
   }

\date{}
 \maketitle

\begin{abstract}
 In solving the problem of asymptotic tracking and disturbance rejection, it has
been long always assumed that the reference to be tracked  and the disturbance
to be rejected must be generated by  an exosystem such as a finite dimensional exosystem with pure
imaginary eigenvalues. The objective of this paper is to solve such a tracking
problem for the unstable Burgers' equation without this assumption. Our treatment of this
problem is straightforward.  Using the method of
   variable   transform,   the tracking problem is split into two separate
   problems: a simple Neumann boundary stabilization problem and a dynamical Neumann
   boundary regulator problem. Unlike the existing literature where the regulator problem is always kept
    independent,   the stabilization problem here is simplified to an independent linear diffusion equation by moving the instability term and the nonlinear term   to the dynamical regulator problem, whereas the dynamical regulator problem does depend on the stabilization problem.
  Thus we can first easily handle the stabilization problem and then
solve the dynamical regulator problem   by using the fundamental theory of
partial differential equations. The boundary feedforward controller is explicitly constructed by using the reference, the disturbance   and the solution of the stabilization problem while the
boundary feedback controller is easily designed for the linear diffusion equation without using a complex method such as the backstepping method. It is proved that, under the designed feedback and feedforward controllers, the tracking error  converges to zero exponentially. This
 theoretical result is confirmed by a numerical example.
\end{abstract}




\section{    Introduction}

Asymptotic tracking and
 disturbance rejection is one of fundamental problems in control theory.
In solving this problem,  it has
been long always assumed that the reference to be tracked  and the disturbance to be
rejected must be generated by an exosystem such as a finite
dimensional exosystem with pure imaginary eigenvalues (see, e.g.,
\cite{Aulisa-book-2015,
 Deutscher2015, Deutscher2016, Deutscher2017, Freudenthaler-2020, Huang-2004-book, Krstic-2009, Liu-book-2010,
 Liu-2021-JMAA, Liu-2020-automatica, Liu-2020-TAC, Liu-2019-MBE, Malchow-2011,
  Meurer-2005, Meurer-2006, Meurer-2011, Meurer-2013, Meurer-2016, Natarajan-2014,
  Utz-2012,Zhou-Weiss2017}). In fact, this assumption is sufficient, but not
  necessary for  making the problem solvable.

The objective of this paper is to solve the tracking problem for the unstable
  Burgers' equation without this assumption:
\begin{eqnarray}\label{bg1}
   \frac{\partial u}{\partial t} &=&\nu   \frac{\partial^2 u}{\partial x^2} -u\frac{\partial u}{\partial x}+au+u_d, \\
 \nu \frac{\partial u}{\partial x}(0,t) &=&f_0,\quad
\nu \frac{\partial u}{\partial x}(1,t)=f_1 , \label{bg2}\\
  u(x,0)&=&u_0(x)  .\label{bg3}
\end{eqnarray}
In the above equations,  $\nu >0$ is a viscosity parameter,
$a= a(x,t)$ is a continuous function,  $u_d= u_d(x,t)$ is a source disturbance,
 $f_0$ and $f_1$
are  control inputs, and the function
 $u_0=u_0(x)$ is an initial state in an appropriate function
space. When $a$ is positively large enough, the equilibrium $0$ of the uncontrolled and undisturbed
Burgers' equation is unstable.

Let $r( t)$ be    a  desired
 reference.  We introduce the   tracking  error
\begin{equation}
  e (  t) =  \int_0^1 u(x, t)dx  -  r(  t). \label{tracking-error}
\end{equation}
 Then the problem of exponential tracking and disturbance rejection for the
   Burgers' equation is to design a feedback and feedforward controllers  $f_0, f_1$
 such that
\begin{equation}\label{error convergence}
   |e ( t)|\le C e^{-\lambda t}
\end{equation}
where $C$ and $\lambda $ are positive constants.

  Mathematical theory on the tracking problem has been well developed for finite
  dimensional control systems (see, e.g., \cite{Huang-2004-book}) and recently extended to partial differential equations (PDEs).
Aulisa \textit{et al.} \cite{Aulisa-book-2015}, Byrnes   \textit{et al.} \cite{Byrnes-2000}, and Natarajan \textit{et al.}
\cite{Natarajan-2014} developed an  abstract theory on the problem in Hilbert spaces and
  applied it to partial differential
equations (PDEs) such as the heat equation and the wave equation.
The flatness method was developed to
handle the trajectory planning and feedforward control design
(see,  e.g., \cite{Freudenthaler-2020, Meurer-2005, Meurer-2006, Meurer-2011, Meurer-2013, Pisano-2011,
 Utz-2012,  Wagner-2009, Zhou-Weiss2017}).
  Wagner \textit{et al.} \cite{Wagner-2009} studied the tracking problem
  for the one-dimensional semilinear wave equation by using
  the flatness method.
    The asymptotic tracking and disturbance rejection of the one-dimensional
parabolic partial differential equations were studied by Deutscher \cite{Deutscher2015, Deutscher2016} and the case with
a long time delay was investigated by Gu \textit{et al}. \cite{Gu2018}. The same problem for the one-dimensional hyperbolic partial differential equations and Schr\"{o}dinger equation were studied by Deutscher \cite{Deutscher2017} and Zhou et al \cite{Zhou-Weiss2017}, respectively. In all of these important researches,
the references and disturbances were assumed to be governed by a finite
dimensional exosystem with pure imaginary eigenvalues.

Although the problem of  feedback stabilization and tracking for the  Burgers' equation
has     received extensive attention (see, e.g., \cite{bk, bka, bgs, ctmk, ik, iy, kr, Krstic-2009, lk00, lk, lmt}), to my knowledge, the tracking problem has not been studied yet in the case of general references and
disturbances, which are not required to be generated by an exosystem.
Our treatment of this
problem is straightforward.
Using the method of     variable   transform used in \cite{Ashley-Liu,
Huang-2004-book, Liu-book-2010, Liu-2021-JMAA, Liu-2020-automatica,
Liu-2020-TAC, Liu-2019-MBE},   the tracking
problem is split into two separate problems:   a simple Neumann boundary
stabilization problem and  a dynamical Neumann boundary regulator problem.
  Unlike the existing literature  where the regulator problem is always kept
    independent (see, e.g., \cite{Aulisa-book-2015, Huang-2004-book}),   the stabilization problem here is simplified to an independent linear diffusion equation by moving the instability term and the nonlinear term   to the dynamical regulator problem, whereas the dynamical regulator problem does depend on the stabilization problem.
  Thus we can first easily handle the stabilization problem and then solve
the dynamical regulator problem   by using the fundamental theory of
partial differential equations (see, e.g., \cite{lsu}).
The feedforward controller is explicitly constructed by using the reference $r$,
 the disturbance $u_d$ and the solution of the stabilization problem while the
 feedback controller is easily designed for the linear diffusion equation without using a complex method such as the backstepping method.
  It is proved that, under the designed feedback and
 feedforward controllers, the tracking error  converges to zero exponentially.
 This theoretical result is confirmed by a numerical example.

\section{Exponential tracking   }

In what follows,
 $H^s(0,1)$   denotes the usual Sobolev space (see \cite{Ada, LM})
 for any $s\in  \mathbb{R}$. For $s\geq 0$, $H^s_0(0,1)$   denotes the completion of
$C_0^\infty(0,1)$ in $H^s(0,1)$, where $C_0^\infty(0,1)$ denotes the space of all
infinitely differentiable functions on $(0,1)$ with compact support  in $(0,1)$.
We use the following $H^1$ norm of $H^1(0,1)$
$$
\|u\|_{H^1}=\left[u(0)^2+\int_0^1 \left(\frac{\partial u}{\partial x}\right)^2dx\right]^{1/2}, \; u\in H^1(0,1),
$$
which is equivalent to the usual one. The  norm on  $L^2(0,1)$ is denoted by $\|\cdot\|$.
It is easy to see that
\begin{equation} \label{l2-h1}
\|u\|^2\le 2 \|u\|_{H^1}^2.
\end{equation}
Let $X$ be a Banach space and $T>0$.
 We denote by $C^n([0,T]; X)$ the  space of $n$ times continuously
differentiable functions defined on $[0,T]$ with values in  $X$, and write
$C([0,T]; X)$ for $C^0([0,T]; X)$.
In what follows, for simplicity,
we omit the indication of the varying range of $x$ and $t$ in equations and we understand that
$x$ varies from $0$ to $1$ and $t$ from $0$ to $\infty$.

To split the tracking problem into   a stabilization
 problem and a dynamical regulator problem, we  introduce the   variable   transform
\begin{equation}\label{linear-transform}
  u =\hat{u} +U,
   \quad f_0   = \hat{f}_0     +F_0,  \quad f_1   = \hat{f}_1     +F_1.
\end{equation}
Subtracting this transform into the problem  \eqref{bg1} - \eqref{bg3} and
the tracking error equation \eqref{tracking-error},
we obtain
\begin{eqnarray*}
\frac{\partial \hat{u}}{\partial t} +\frac{\partial U}{\partial t}   &=&
\nu\frac{\partial^2 \hat{u}}{\partial x^2} +
\nu\frac{\partial^2 U}{\partial x^2} +a(\hat{u}+U) \\
&&- (\hat{u} +U)\left(\frac{\partial  \hat{u}}{\partial x } +
 \frac{\partial  U}{\partial x }\right) +u_d , \\
\nu \frac{\partial \hat{u}}{\partial x} +\nu \frac{\partial U}{\partial
 x}(0,t) &=&\hat{f}_0   +F_0 ,\\
 \nu \frac{\partial \hat{u}}{\partial x} +\nu \frac{\partial U}{\partial
 x}(1,t) &=&\hat{f}_1   +F_1 ,\\
 \hat{u}(x,0)+U(x,0)
  &=&u_0(x) ,  \\
  e(t)&=& \int_0^1  [\hat{u}(x, t)+U(x,t)]dx
  - r(t).
\end{eqnarray*}
This problem can be split into
 a stabilization problem
\begin{eqnarray}
\frac{\partial \hat{u}}{\partial t}     &=&
\nu\frac{\partial^2 \hat{u}}{\partial x^2} -
   \hat{u}  \frac{\partial  \hat{u}}{\partial x }   ,\label{transformed-bg-eq} \\
\nu \frac{\partial \hat{u}}{\partial x} (0,t)  &=&\hat{f}_0    ,\quad
 \nu \frac{\partial \hat{u}}{\partial x} (1,t)  =\hat{f}_1    ,\label{transformed-BC-eq}\\
  \hat{u}(x,0)
  &=&u_0(x)- r(0),  \label{transformed-IC-eq}\\
  e(t)&=&  \int_0^1   \hat{u}(x, t) dx,
\label{transformed-error-eq}
\end{eqnarray}
and a dynamical regulator problem
\begin{eqnarray}
 \frac{\partial U}{\partial t}   &=&
\nu\frac{\partial^2 U}{\partial x^2}  -  U
 \frac{\partial  U}{\partial x }+a(\hat{u}+U)  \nonumber\\
 &&-  \left(\hat{u}
 \frac{\partial  U}{\partial x }  +  U  \frac{\partial  \hat{u}}{\partial x }\right)+u_d   , \label{regulator1}\\
 \nu  \frac{\partial U}{\partial
 x}(0,t) &=& F_0 ,\quad
  \nu \frac{\partial U}{\partial
 x}(1,t)  = F_1
  \label{regulator-bc}\\
    U(x,0)
  &=& r(0),  \label{regulator-ic}\\
  \int_0^1  U(x,t) dx &=& r(t).  \label{regulator2}
\end{eqnarray}

 In the control design (see, e.g., \cite{Aulisa-book-2015, Huang-2004-book}),
 the terms $a(\hat{u}+U)$ and  $\hat{u}
 \frac{\partial  U}{\partial x }  +  U  \frac{\partial  \hat{u}}{\partial x }$
  are usually put in the equation \eqref{transformed-bg-eq} in the stabilization problem such that the regulator problem is independent from the stabilization
 problem. However, if we do so for the Burgers' equation, then the feedback controllers $\hat{f}_0  $ and $\hat{f}_1  $
for the stabilization problem \eqref{transformed-bg-eq} - \eqref{transformed-IC-eq}
is difficult or impossible to design. In fact, we can further simplify the stabilization problem by moving the term $\hat{u}  \frac{\partial  \hat{u}}{\partial x }$ from the stabilization problem to the regulator problem and then obtain the following linear stabilization problem and the regulator problem:
\begin{eqnarray}
\frac{\partial \hat{u}}{\partial t}     &=&
\nu\frac{\partial^2 \hat{u}}{\partial x^2}     ,\label{transformed-bg-eq-new} \\
\nu \frac{\partial \hat{u}}{\partial x} (0,t)  &=&\hat{f}_0    ,\quad
 \nu \frac{\partial \hat{u}}{\partial x} (1,t)  =\hat{f}_1    ,\label{transformed-BC-eq-new}\\
  \hat{u}(x,0)
  &=&u_0(x)- r(0),  \label{transformed-IC-eq-new}\\
  e(t)&=&  \int_0^1   \hat{u}(x, t) dx,
\label{transformed-error-eq-new}
\end{eqnarray}
and
\begin{eqnarray}
 \frac{\partial U}{\partial t}   &=&
\nu\frac{\partial^2 U}{\partial x^2}  -  U
 \frac{\partial  U}{\partial x }+a(\hat{u}+U)  \nonumber\\
 &&-
   \hat{u}  \frac{\partial  \hat{u}}{\partial x }-  \left(\hat{u}
 \frac{\partial  U}{\partial x }  +  U  \frac{\partial  \hat{u}}{\partial x }\right)+u_d   , \label{regulator1-new}\\
 \nu  \frac{\partial U}{\partial
 x}(0,t) &=& F_0 ,\quad
  \nu \frac{\partial U}{\partial
 x}(1,t)  = F_1
  \label{regulator-bc-new}\\
    U(x,0)
  &=& r(0),  \label{regulator-ic-new}\\
  \int_0^1  U(x,t) dx &=& r(t).  \label{regulator2-new}
\end{eqnarray}

Because the regulator problem \eqref{regulator1} -
\eqref{regulator2} of \eqref{regulator1-new} -
\eqref{regulator2-new} is not dissipative, its well-posedness
is challenging and open.  Since $\int_0^1  U(x,t) dx = r(t)$
   exists   for all times, we could conjecture that it has a unique global solution. We will use this conjecture in the
   following theorem.

Using the stabilization problem \eqref{transformed-bg-eq} - \eqref{transformed-IC-eq}
 and the regulator problem \eqref{regulator1} - \eqref{regulator2}, we can design the feedback
 and feedforward controllers as stated in the following theorem.

\begin{theorem} \label{main-theorem1} Assume that $k>1/6$  and the
initial condition $u_0\in H^2(0,1)$.  Suppose that
  $a(x,t)$ and $u_d(x, t)$ are continuous and   $r(  t)$ is continuously differentiable.
  Let $\hat{u}$ be the solution of the stabilization problem \eqref{transformed-bg-eq} - \eqref{transformed-IC-eq} and assume that the regulator problem \eqref{regulator1} - \eqref{regulator2} has a unique global classical solution.
   Then, under the feedback and feedforward controllers:
   \begin{eqnarray}
   \hat{f}_0(t)  &=&  k\left[\hat{u} (0,t)+[\hat{u} (0,t) ]^3\right]  ,
    \label{controller-hf0} \\
    \hat{f}_1(t)  &=&-k\left[\hat{u}(1,t) +[\hat{u}(1,t) ]^3\right],
    \label{controller-hf1}\\
F_0(t)  &=&    \hat{u} (0,t) U (0,t)  +\frac{1}{2}[U(0,t)]^2 ,
    \label{controller-F0} \\
    F_1(t)  &=&  \hat{u}(1,t) U(1,t)  +\frac{1}{2}[U(1,t)]^2 +r'(t) \nonumber\\
 &&- \int_0^1[a(x,t)(\hat{u}(x,t)+U(x,t))+ u_d(x,t)]dx .
    \label{controller-F1}
 \end{eqnarray}
  the problem \eqref{bg1} - \eqref{bg3} has a unique solution satisfying
  \begin{equation}\label{exp-tracking}
   | e( t)|
   \le \|u_0-r(0)\|e^{-\lambda t/2},
  \end{equation}
where  $\lambda =\min\left(\nu, k-\frac{1}{6}\right)$.
\end{theorem}

\emph{Proof}.   If  $k>1/6$  and the
initial condition $u_0\in H^2(0,1)$, then it was proved in  \cite{kr, lk00} that
  the stabilization problem \eqref{transformed-bg-eq} - \eqref{transformed-IC-eq} with the feedback controllers \eqref{controller-hf0}
and \eqref{controller-hf1} has a unique solution
satisfying
$$\hat{u}\in C([0,\infty); H^2(0,1)).$$
Moreover, multiplying the equation \eqref{transformed-bg-eq} by
$\hat{u}$ and integrating it from $0$ to $1$, we obtain
\begin{eqnarray*}
 && \frac{1}{2}\frac{d}{dt}\int_0^1[\hat{u}(x,t)]^2dx \\
 &=&
    -k\left[\hat{u}(1,t) +[\hat{u}(1,t) ]^3\right]\hat{u}(1,t)\\
    &&-k\left[\hat{u} (0,t)+[\hat{u} (0,t) ]^3\right]\hat{u} (0,t)\\
    &&-\nu\int_0^1\left[\frac{\partial \hat{u}}{\partial x}(x,t)\right]^2dx
    -\frac{1}{3}\left([\hat{u}(1,t) ]^3-[\hat{u}(0,t) ]^3\right)\\
    &\le& -\left(k-\frac{1}{6}\right)\left([\hat{u}(0,t)]^2 +[\hat{u}(1,t)]^2
    \right.\\
    &&\left. +[\hat{u}(0,t)]^4 +[\hat{u}(1,t) ]^4\right)-\nu\int_0^1\left[\frac{\partial \hat{u}}{\partial x}(x,t)\right]^2dx\\
    &\le& -\left(k-\frac{1}{6}\right) [\hat{u}(0,t)]^2  -\nu\int_0^1\left[\frac{\partial \hat{u}}{\partial x}(x,t)\right]^2dx\\
    &\le& -\lambda\left([\hat{u}(0,t)]^2   +\int_0^1\left[\frac{\partial \hat{u}}{\partial x}(x,t)\right]^2dx\right)\\
    &\le& -\frac{\lambda}{2}  \int_0^1\left[  \hat{u} (x,t)\right]^2dx.
    \quad (\mbox{use } \eqref{l2-h1})
\end{eqnarray*}
Solving this inequality, we obtain
  \begin{equation}\label{exp-decay}
   \| \hat{u}( t)\|
   \le \|u_0- r(0)\| e^{-\lambda t/2}.
  \end{equation}


Integrating the  equation  \eqref{regulator1} over $[0,1]$ , we obtain
\begin{eqnarray*}
 && \frac{d}{dt}\int_0^1 U(x,t)dx = \nu \frac{\partial U}{\partial
 x}(1,t)-\nu \frac{\partial U}{\partial
 x}(0,t)\\
 && +
  \hat{u}(0,t)U(0,t)+\frac{1}{2}U^2(0,t)\\
   &&- \hat{u}(1,t)U(1,t)-\frac{1}{2}U^2(1,t)+\int_0^1 u_d(x, t) dx \\
   &&+ \int_0^1 a(x,t)[\hat{u}(x,t)+U(x,t) ]dx.
\end{eqnarray*}
It then follows from the boundary condition \eqref{regulator-bc}
and  the   feedforward controllers \eqref{controller-F0}
and \eqref{controller-F1}
  that
\begin{equation*}
  \frac{d}{dt}\int_0^1 U(x,t)dx =  r'(t)  .
\end{equation*}
Integrating this equation and using the initial condition
\eqref{regulator-ic}, we obtain
\begin{equation}
\int_0^1  U(x,t) dx = r(t).
\end{equation}
So $U$ satisfies the equation \eqref{regulator2}.
Since $ \int_0^1 U(x, t)dx$ exists
for all $t\ge 0$,
the problem \eqref{regulator1} - \eqref{regulator-ic}
has a unique solution $U$   for all $t>0$.

Finally it     follows from  the equations   \eqref{transformed-error-eq},
\eqref{regulator2},  and \eqref{exp-decay} that
\begin{eqnarray*}
|e(t)| &=& \left| \int_0^1 u(x, t)dx -r(t)\right|\nonumber\\
  &=& \left| \int_0^1 u(x, t)dx -\int_0^1  U(x,t) dx\right|\nonumber\\
 &=& \left| \int_0^1\hat{u}(x, t)dx \right|\nonumber\\
&\le & \left[\int_0^1 dx \right]^{1/2}  \left[\int_0^1  [ \hat{u}(x,  t)]^2dx \right]^{1/2}\nonumber\\
 &\le &\|u_0- r(0)\|e^{-\lambda t /2}.
\end{eqnarray*}
This completes the proof.

Using the stabilization problem \eqref{transformed-bg-eq-new} - \eqref{transformed-IC-eq-new}
 and the regulator problem \eqref{regulator1-new} - \eqref{regulator2-new}, we can design the feedback
 and feedforward controllers as stated in the following theorem.

\begin{theorem} \label{main-theorem2} Assume that $k>0$  and the
initial condition $u_0\in H^2(0,1)$.  Suppose that
  $a(x,t)$ and $u_d(x, t)$ are continuous and   $r(  t)$ is continuously differentiable.
  Let $\hat{u}$ be the solution of the stabilization problem \eqref{transformed-bg-eq-new} - \eqref{transformed-IC-eq-new} and assume that the regulator problem \eqref{regulator1-new} - \eqref{regulator2-new} has a unique global classical solution.
   Then, under the feedback and feedforward controllers:
   \begin{eqnarray}
   \hat{f}_0(t)  &=&  k \hat{u} (0,t)   ,
    \label{controller-hf0-new} \\
    \hat{f}_1(t)  &=&-k \hat{u}(1,t)  ,
    \label{controller-hf1-new}\\
F_0(t)  &=&    \hat{u} (0,t) U (0,t)  +\frac{1}{2}\left([\hat{u} (0,t) ]^2+[U(0,t)]^2\right) ,
    \label{controller-F0-new} \\
    F_1(t)  &=&  \hat{u}(1,t) U(1,t)  +\frac{1}{2}\left([\hat{u} (1,t) ]^2+[U(1,t)]^2 \right)+r'(t) \nonumber\\
 &&- \int_0^1[a(x,t)(\hat{u}(x,t)+U(x,t))+ u_d(x,t)]dx .
    \label{controller-F1-new}
 \end{eqnarray}
  the problem \eqref{bg1} - \eqref{bg3} has a unique solution satisfying
  \begin{equation}\label{exp-tracking}
   | e( t)|
   \le \|u_0-r(0)\|e^{-\lambda t/2},
  \end{equation}
where  $\lambda =\min\left(\nu, k \right)$.
\end{theorem}

The proof of this theorem is the same as the proof of Theorem \ref{main-theorem1}.

\section{A numerical example}

\begin{figure}[t]
 \begin{center}
\includegraphics[width=4.5cm, height=4cm] 
{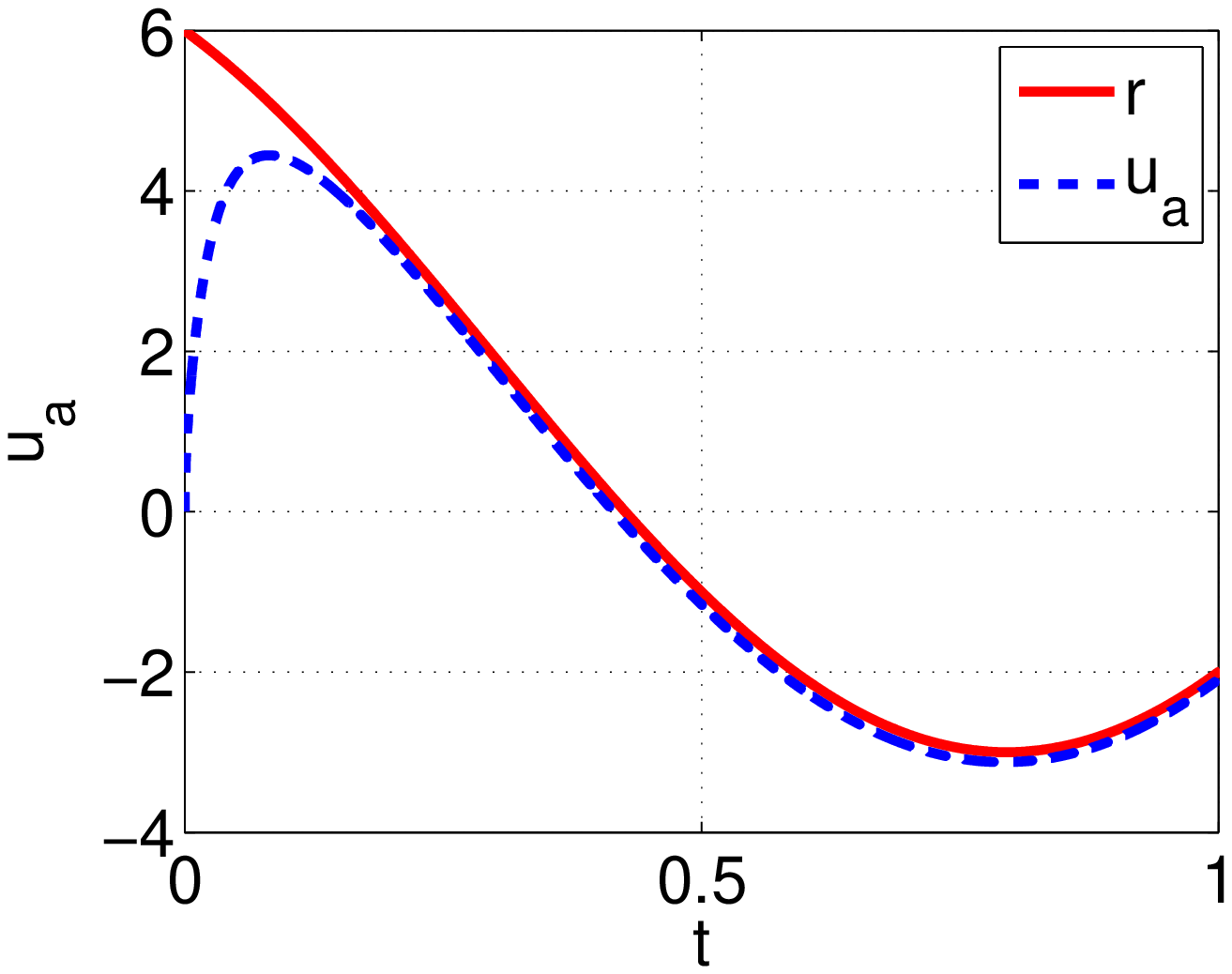}
\includegraphics[width=4.5cm, height=4cm] 
{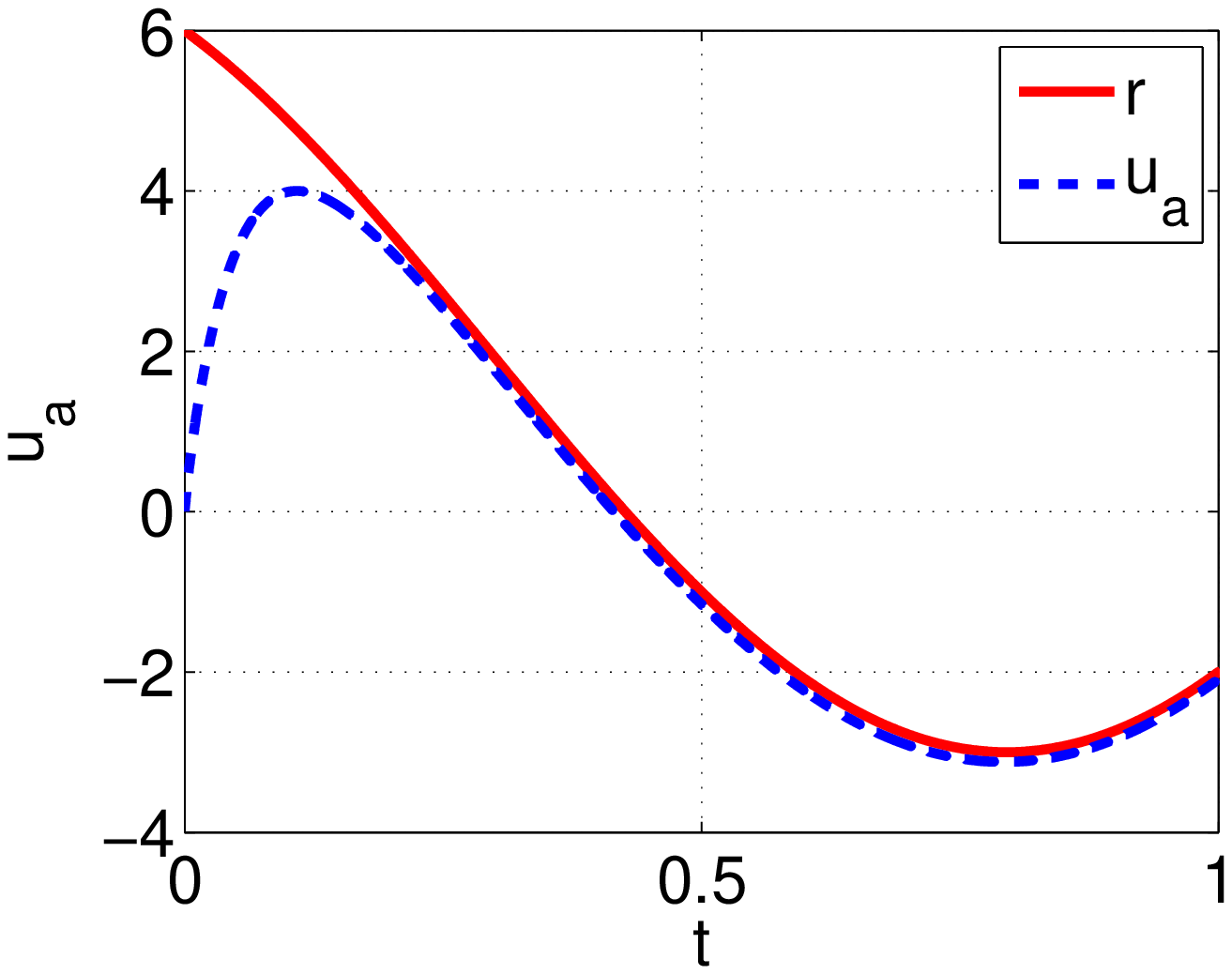}
 \caption{The average of $u$, $u_a(t) =\int_0^1 u(x,t)dx$,      quickly tracks the     reference $r(t) = 2+4\cos(\pi t)-3\sin(\pi t)$ under   the feedback and feedforward controllers either \eqref{controller-hf0} - \eqref{controller-F1} (left figure) or \eqref{controller-hf0-new} - \eqref{controller-F1-new} (right figure).}
 \label{average-tracking-fig}
 \end{center}
\end{figure}

  We conduct a numerical simulation to confirm the above
  theoretical result.
  In the numerical computations, we
  take      $\nu = 5$, $a=20$, $ k= 15$, $u_0(x) =0$,
   $r(t) = 2+4\cos(\pi t)-3\sin(\pi t)$, and
   $u_d(x, t) = 3+5\cos(\pi x)\sin(\pi t)-
   2\sin(\pi x)\cos(\pi t)$.
Then the problem \eqref{bg1} - \eqref{bg3} and the problem
\eqref{transformed-bg-eq} - \eqref{transformed-IC-eq}  are solved numerically by
the difference method.
  The   Figure \ref{average-tracking-fig} shows that the average of $u$, $u_a(t) =\int_0^1 u(x,t)dx$,      quickly tracks the     reference $r(t) = 2+\cos(\pi t)-3\sin(\pi t)$ under   the feedback and feedforward controllers either \eqref{controller-hf0} - \eqref{controller-F1} or \eqref{controller-hf0-new} - \eqref{controller-F1-new}.

\section{Discussion}

Other tracking errors can be considered.
Let $r_0(t)$ and $r_1(t)$ be two desired references and define the tracking
error by
\begin{equation}\label{tracking-error1}
  e_1(t) = u(0,t)-r_0(t) +u(1, t)-r_1(t).
\end{equation}
Using  the variable transform \eqref{linear-transform} and this tracking error, we can obtain
 the stabilization problem
\begin{eqnarray}
\frac{\partial \hat{u}}{\partial t}     &=&
\nu\frac{\partial^2 \hat{u}}{\partial x^2} -
   \hat{u}  \frac{\partial  \hat{u}}{\partial x }   ,\label{error1-transformed-bg-eq} \\
\nu\frac{\partial \hat{u}}{\partial x} (0,t)  &=&\hat{f}_0    ,\quad
 \nu\frac{\partial \hat{u}}{\partial x} (1,t)  =\hat{f}_1    ,\label{error1-transformed-BC-eq}\\
  \hat{u}(x,0)
  &=&u_0(x)-r_0(0)(1-x)-r_1(0)x,  \label{error1-transformed-IC-eq}\\
  e_1(t)&=&    \hat{u}(0, t) +\hat{u}(1, t),
\label{error1-transformed-error-eq}
\end{eqnarray}
and the dynamical regulator problem
\begin{eqnarray}
 \frac{\partial U}{\partial t}   &=&
\nu\frac{\partial^2 U}{\partial x^2}   -  U
 \frac{\partial  U}{\partial x }+a(\hat{u}+U)  \nonumber\\
 &&-  \left(\hat{u}
 \frac{\partial  U}{\partial x }  +  U  \frac{\partial  \hat{u}}{\partial x }\right)+u_d   , \label{error1-regulator1}\\
  U(0,t)  &=& r_0(t), \quad  U(1,t)   =r_1(t), \label{error1-regulator2}\\
    U(x,0)
  &=& r_0(0)(1-x)+r_1(0)x,  \label{error1-regulator-ic}\\
 \nu \frac{\partial U}{\partial
 x}(0,t) &=& F_0 ,\quad
 \nu \frac{\partial U}{\partial
 x}(1,t) =F_1 ,
  \label{error1-regulator-bc}
\end{eqnarray}
If the problem \eqref{error1-regulator1} - \eqref{error1-regulator-ic} has
a global solution, then the feedforward controllers are given by the
equation \eqref{error1-regulator-bc}. However, it seems challenging to show it
has a global solution even though it has a unique classical solution within some time from $0$ to $T$ (see, e.g., \cite{lsu}).

We can also consider   the tracking
error
\begin{equation}\label{tracking-error2}
  e_2(t) = u(0,t)-r_0(t) +\int_0^1u(x,t)dx-r(t).
\end{equation}
Using  the variable transform \eqref{linear-transform}   and this tracking error, we can obtain
 the stabilization problem
\begin{eqnarray}
\frac{\partial \hat{u}}{\partial t}     &=&
\nu\frac{\partial^2 \hat{u}}{\partial x^2} -
   \hat{u}  \frac{\partial  \hat{u}}{\partial x }   ,\label{error2-transformed-bg-eq} \\
\nu\frac{\partial \hat{u}}{\partial x} (0,t)  &=&\hat{f}_0    ,\quad
 \nu\frac{\partial \hat{u}}{\partial x} (1,t)  =\hat{f}_1    ,\label{error2-transformed-BC-eq}\\
  \hat{u}(x,0)
  &=&u_0(x)-u_r(0),  \label{error2-transformed-IC-eq}\\
  e_2(t)&=&     \hat{u}(0, t) +\int_0^1   \hat{u}(x, t) dx,
\label{error2-transformed-error-eq}
\end{eqnarray}
and the dynamical regulator problem
\begin{eqnarray}
 \frac{\partial U}{\partial t}   &=&
\nu\frac{\partial^2 U}{\partial x^2}  -  U
 \frac{\partial  U}{\partial x }+a(\hat{u}+U)  \nonumber\\
 &&-  \left(\hat{u}
 \frac{\partial  U}{\partial x }  +  U  \frac{\partial  \hat{u}}{\partial x }\right)+u_d   , \label{error2-regulator1}\\
   U(0,t)   &=& r_0(t)  ,\quad
 \nu \frac{\partial U}{\partial
 x}(1,t) =F_1 ,
  \label{error2-regulator-bc}\\
    U(x,0)
  &=& r(0),  \label{error2-regulator-ic}\\
      \int_0^1  U(x,t) dx &=& r(t),
       \label{error2-regulator2}\\
      \nu \frac{\partial U}{\partial
 x}(0,t) &=& F_0. \label{error2-controller-F0}
\end{eqnarray}
To find $F_1$, we integrate the equation \eqref{error2-regulator1}
from $0$ to $1$ and use the equation \eqref{error2-regulator2} to obtain
\begin{eqnarray}
  F_1 &=&  \nu \frac{\partial U}{\partial
 x}(0,t)+\frac{1}{2}[(U(1,t))^2-(U(0,t))^2]\nonumber\\
 &&+\hat{u}(1,t)U(1,t)-\hat{u}(0,t)U(0,t) +r'(t) \nonumber  \\
  && - \int_0^1[a(x,t)(\hat{u}(x,t)+U(x,t)+ u_d(x,t)]dx.
\end{eqnarray}
This results the following complex boundary value problem
\begin{eqnarray}
 \frac{\partial U}{\partial t}   &=&
\nu\frac{\partial^2 U}{\partial x^2}  - (\hat{u} +U)
 \frac{\partial  U}{\partial x }  -  U  \frac{\partial  \hat{u}}{\partial x }+u_d   , \label{openregulator1}\\
   U(0,t)   &=& r_0(t)  ,\\
 \nu \frac{\partial U}{\partial
 x}(1,t) &=&\nu \frac{\partial U}{\partial
 x}(0,t)+\frac{1}{2}[(U(1,t))^2-(U(0,t))^2]\nonumber\\
 &&+\hat{u}(1,t)U(1,t)-\hat{u}(0,t)U(0,t) +r'(t)  \nonumber  \\
  && -\int_0^1[a(x,t)(\hat{u}(x,t)+U(x,t)+ u_d(x,t)]dx ,
  \label{open-regulator-bc}\\
    U(x,0)
  &=& r(0).  \label{open-regulator-ic}
\end{eqnarray}
It seems that the proof of the solution existence of the problem is challenging.

\end{document}